\documentclass[conference]{IEEEtran}

\usepackage{cite}
\usepackage{amsmath,amssymb,amsfonts}
\usepackage{algorithmic}
\usepackage{mathtools}

\usepackage{textcomp}
\def\BibTeX{{\rm B\kern-.05em{\sc i\kern-.025em b}\kern-.08em
    T\kern-.1667em\lower.7ex\hbox{E}\kern-.125emX}}

\IEEEoverridecommandlockouts
\usepackage{graphicx}
\usepackage{cite}
\usepackage{picinpar}
\usepackage{amsmath}
\usepackage{url}
\usepackage{flushend}
\usepackage[latin1]{inputenc}
\usepackage{colortbl}
\usepackage{soul}
\usepackage{multirow}
\usepackage{pifont}
\usepackage{color}
\usepackage{alltt}
\usepackage[hidelinks]{hyperref}
\usepackage{enumerate}
\usepackage{siunitx}

\usepackage{epstopdf}
\usepackage{pbox}
\usepackage{verbatim}
\usepackage{leftidx}
\usepackage{amssymb}
\usepackage{xcolor}
\usepackage{listings}
\usepackage{textcomp}

\newtheorem{Theorem}{Theorem}[section]

\newtheorem{Corollary}{Corollary}[section]
\newtheorem{Lemma}[Theorem]{Lemma}
\newtheorem{Definition}{Definition}[section]

 \newcommand{\R}{\mathbb{R}}

  \newcommand{\zero}{0}

\newcommand{\RNum}[1]{\uppercase\expandafter{\romannumeral #1\relax}}
\graphicspath{ {Pictures/} }

\def\BibTeX{{\rm B\kern-.05em{\sc i\kern-.025em b}\kern-.08em
    T\kern-.1667em\lower.7ex\hbox{E}\kern-.125emX}}
    
\usepackage{comment}

\begin{document}
\title{Optimal Rejection of Bounded Perturbations in Linear Leader-Following Consensus Protocol: Method Invariant Ellipsoid} 
\author{Siyuan WANG,   Andrey POLYAKOV, Min LI, Gang ZHENG,  Driss BOUTAT

\thanks{This work is supported by the Project FEDER Medibot, OpMedibot and  the China Scholarship Council (CSC) under Grant 202106160022. }
\thanks{Siyuan WANG and Driss Boutat  are with INSA Centre Val de Loire, Bourges, France (e-mail: Siyuan.wang@insa-cvl.fr, Driss.Boutat@insa-cvl.fr). }
\thanks{Andrey POLYAKOV and Gang ZHENG are with INRIA Lille, Lille, France. (e-mail: andrey.polyakov@inria.fr, gang.zheng@inria.fr).}
\thanks{Min Li is the corresponding author with 
INRIA Lille, Lille, France (e-mail: min.li@inria.fr).}
}


\maketitle

\begin{abstract}

The objective of the invariant ellipsoid method is to minimize the smallest invariant and attractive set of a linear control system operating under the influence of bounded external disturbances. In this paper,  this method is extended into the leader-following consensus problem. Initially,  a linear control protocol is designed for the Multi-agent System  without  disturbances. Subsequently, in the presence of bounded disturbances, by employing a similar linear control protocol,   a necessary and sufficient condition is introduced to derive the optimal control parameters for the MAS such that the state of followers converge and remain in an minimal invariant ellipsoid around the state of the leader.

\end{abstract}


\markboth{IEEE Transactions }%
{}

\definecolor{limegreen}{rgb}{0.2, 0.8, 0.2}
\definecolor{forestgreen}{rgb}{0.13, 0.55, 0.13}
\definecolor{greenhtml}{rgb}{0.0, 0.5, 0.0}

\section{Introduction} \label{sec:intro}

The Multi-Agent System (MAS) is a kind of system consisted of multiple agents, which could achieve some global behavior via local communication. During task execution, MAS performs better working
efficiency, reduced sensitivity and increased flexibility. Thus MAS  has received considerable academic attention, and has be researched actively on control problem such as consensus \cite{olfati2004consensus}, formation \cite{mastellone2008formation,huo2021pigeon} and containment \cite{li2013distributed}; as well on estimation problem such as distributed estimation \cite{park2016design}. The obtained theory has been considered into some application scenarios such as swarm of mobile robots \cite{9976029}, network security \cite{9285180} and smart grid \cite{9076711}. 

Consensus problem is one of the most fundamental research issue with MAS, which can be subdivided into two categories, the
leaderless consensus and the leader-following consensus. In the case of leaderless consensus, the ultimate shared  position can be not pre-selected. However, in certain scenarios, as demonstrated in  \cite{ni2010leader},  the leader-following problem might be more useful wherein a virtual or real leader can be introduced to guide the MAS in following a predefined or unknown trajectory. 



From a practical application perspective, while formulating MAS sonsensus problems, discrepancies inevitably arise between the actual system and the mathematical model employed for control design. These discrepancies stem from unmodelled dynamics, uncertainties in MAS parameters, or the simplification of complex plant dynamics. Therefore, how to guarantee the required performance of MAS under certain (matched and mismatched) disturbances is one of the important issues for the practical applications. This challenge has led to the development of robust control methods aimed at addressing this issue.

Robust control design of MAS is an approach focused on mitigating disturbances within a nominal system. Its main objective is to attain a specified level of performance or stability in the system, even in the presence of bounded disturbances. Numerous established methodologies for robust control design have been developed, including techniques such as sliding mode control \cite{utkin1977variable,waslander2005multi}, $H_\infty$ approach \cite{kimura1984robust,li2011h} , attractive/invariant ellipsoid method \cite{nazin2007rejection,deng2016consensus}, among others.

The Sliding Mode Control (SMC) as a robust control design methodology has been recognized since the 1960s in Russia. V. Utkin published the first survey in 1977 \cite{utkin1977variable}. Sliding mode method could provide a good performance especially for the system with bounded matched disturbance, since this kind of disturbance can be compensated by the discontinuity of the SMC \cite{utkin2013sliding}, \cite{shtessel2014sliding}\cite{utkin2017sliding}. However, SMC may not provide a satisfied performance for the system with unmatched disturbance \cite{rubagotti2011integral}. 
$H_\infty$ techniques offer distinct advantages over classical control methods, particularly in their applicability to problems  of  involving multivariable systems
with cross-coupling between channels \cite{zames1981feedback,orlov2014advanced,kimura1984robust,gershon2005h}. However, it's important to acknowledge that the use of $H_\infty$ techniques comes with certain challenges. These challenges include the requirement for a high level of mathematical proficiency to effectively apply these techniques and the necessity for a reasonably accurate model of the system being controlled. 

The basic idea of the attractive/invariant method were introduced in \cite{bertsekas1971recursive,glover1971control,schweppe1968recursive}. This method has found widespread application in addressing a range of control and estimation challenges, both for linear systems \cite{khlebnikov2011optimization} and nonlinear plant models \cite{Poznyak_etal2014:Book}. One  important characteristic of this method is that the system state converges to the minimal ellipsoid regardless of perturbations or uncertainty,  provided they satisfy certain bounds. The utilization of invariant ellipsoids has facilitated the reformulation of original problems in terms of matrix inequalities, thereby reducing optimal control design problems to solving semi-definite programs. \textcolor{blue}{Many researchers have applied the invariant set method to solve the problems in MAS. A notable advancement in the field of MAS and invariant set methodology is presented in the study by Deng et al. \cite{deng2016consensus}. It provides a sufficient condition of being an invariant set for MAS. In \cite{beji2015invariant},  a finite-time control protocol of reaching the invariant set is designed. In \cite{she2020invariant}, the invariant set method is
proposed to analyze the limit behaviors of all the synchronous trajectories. However, these  studies above primarily provides or use  the sufficient condition for the establishment of the invariant set in MAS, which means that the proposed control protocol may not be optimal.}
\textcolor{blue}{
This article's primary contribution lies in expanding the attractive/invariant ellipsoid method to MAS. We propose a necessary and sufficient condition that ensures the control protocol is optimal, specifically designed to minimize the impact of external disturbances. Furthermore, we derive a sufficient condition for the characterization of an attractive ellipsoid.
}

The rest of paper is structured as follows: In Section \ref{Preliminary}, we delve into the concepts of graph theories and invariant set theories. Section \ref{sec:Pro_statement} outlines the problem statement. Section \ref{Main_results} explains the main results regarding the necessary and sufficient conditions for being an attractive/invariant ellipsoid of the MAS under bounded disturbance. Finally, Section \ref{sim_results} presents simulation results that support the proposed theories.

{\itshape Notation :}
$\mathbb{R}$ is the  set of real numbers; 
$\mathbb{R}_+ = \{x\in \mathbb{R}:x> 0\}$;
$\|x\| =\sqrt{x^\top x}  $ is the Euclidean norm in $\R^n$;  
 $\|x\|_P = \sqrt{x^\top P x}$ is the  weight norm in $\mathbb{R}^n$;
$P\succ0 (\prec0,\succeq 0,\preceq 0)$ for a symmetric matrix $P\in \mathbb{R}^{n\times n}$ means that the matrix $P$ is positive (negative) definite (semi-definite);  
$\lambda_{min}(P)$ and $\lambda_{max}(P)$ represent the minimal and maximal eigenvalue  of  matrix $P$, respectively;
 for $P\succeq 0$ the square root of $P$ is a matrix $M=P^{\frac{1}{2}}$ such that  $M^2=P$; $\mathrm{tr}(P)$ denotes the  trace of matrix $P$; 
 $L^{\infty}$ is the space of Lebesgue measurable essentially bounded function $\sigma : \R_+ \mapsto \R^n$ with norm defined as $\|\sigma\|_{L^\infty} : =\mathrm{ess} \sup_{t\in\R_+}\|\sigma(t)\|_\infty<+\infty$; $1_N \in \mathbb{R}^N$ is a vector with all elements equal to 1.

\section{Preliminary}
\label{Preliminary}
\subsection{Graph Theory}
A fixed graph is usually characterized by  $\mathcal{G}=\{\mathcal{V},\mathcal{E},\mathcal{C}\}$, where $\mathcal{V}=\{0,\dots,N\}$ is the node set;
$\mathcal{E}=\{(i,j)|i,j\in \mathcal{V}\}$ is the edge set, $(i,j)\in \mathcal{E}$ if the local information of node $j$ could convey to node $i$, and node $j$ makes a neighbor of node $i$, $n_i$ denotes the quantity of neighbors of node $i$;
$\mathcal{C}$ is the weighted adjacency matrix whose elements are denoted by $c_{ij},\;i,j\in\mathcal{V}$, $c_{ij}>0$ if $(i,j)\in \mathcal{E}$, and $c_{ij}=0$ otherwise. The self-loop phenomenon is excluded in this work, {\it{i.e.,} }$c_{ii}=0$. The graph   $\mathcal{G}$ is directed if $(i,j)\in \mathcal{E}$ does not imply $(j,i)\in \mathcal{E}$, while $\mathcal{G}$ is undirected if $(i,j)\in \mathcal{E}\Leftrightarrow(j,i)\in \mathcal{E}$. A directed path from node $i$ to node $j$ in the graph $\mathcal{G}$ is a sequence of nodes $\overline{i_0, i_s}$, where $i_0=i$, $i_s=j$ and $(i_{\kappa+1},i_\kappa)\in\mathcal{E}$, $\kappa=\overline{0,s-1}$. 
A directed tree is a directed graph with a node called root, which possesses a unique directed path to each of the rest node. $\mathcal{T}_{\mathcal{G}}=\{\mathcal{V}_{\mathcal{T}},\mathcal{E}_{\mathcal{T}},\mathcal{C}_{\mathcal{T}}\}$ is a spanning tree of graph $\mathcal{G}$ if $\mathcal{T}_{\mathcal{G}}$ is a directed tree and $\mathcal{V}_{\mathcal{T}}=\mathcal{V}$.
The Laplacian matrix $\mathcal{L}$, induced by  $\mathcal{C}$, with elements $l_{ij},~i,j\in\mathcal{V}$ defined as
\begin{equation*}
  l_{ij}=\left\{
          \begin{array}{lll}
           -c_{ij}  & \text{if} &\quad i\neq j, \\
              \sum_{k=0}^Nc_{ik}&\text{if} &\quad i=j. \\
           \end{array}
         \right.
\end{equation*}
Thereby it is obvious one Laplacian matrix is corresponding to one specified $\mathcal{G}$. The Laplacian matrix has some particular properties, {\it{e.g.,}} 
all eigenvalues of $\mathcal{L}$ have non-negative real parts; if $\mathcal{G}$ is an undirected graph, then the associated Laplacian matrix $\mathcal{L}\succeq0$. A directed graph $\mathcal{G}$ has a spanning tree if and only if the associated Laplacian matrix $\mathcal{L}$ has exactly one zero eigenvalue \cite{bullo2009distributed}.
\begin{Lemma}[\cite{GUAN20121397}] \label{lem:tildeL}
Let
$$\hat{\mathcal{L}}=\left[\begin{matrix}
l_{11}-l_{01} & \ldots & l_{1 N}-l_{0N}\\
\vdots        & \ddots &  \vdots\\
l_{N}-l_{01}  & \ldots & l_{NN}-l_{0N}  \\
\end{matrix}\right],$$ 
 $\zeta_i$, $i=1,2,\dots,N$ be the eigenvalues of $\hat{\mathcal{L}}$, and $\mu_j$, $j=0,1,\dots,N$ be the eigenvalues of ${\mathcal{L}}$, respectively, where $|\zeta_1|\leq|\zeta_2|\leq\dots\leq|\zeta_N|$ and $0=|\mu_0|\leq|\mu_1|\leq\dots\leq|\mu_N|$. Then, $\zeta_1=\mu_1,\dots,\zeta_{N}=\mu_{N}$.
\end{Lemma}



\subsection{Invariant Set Theory}


 Recall the definition of attractive/invariant ellipsoid.

\begin{Definition}[\cite{khlebnikov2011optimization} ]
\it \label{def:lin_att_inv}
Let symmetric matrix $P\succ 0$,     the ellipsoid 
    \begin{equation}
        \varepsilon(P) = \{e\in \R^{Nn} : \| e\|_P \leq 1\},\quad P\succ  0
    \end{equation}
    centered at origin is called 
    \begin{itemize}
        \item  \textbf{Invariant} for system \eqref{eq:mas_0}-\eqref{eq:mas_N}, if  $e(0)\in \varepsilon(P)$ implies $e(t)\in\varepsilon(P)$ for $\forall t\geq 0$.
        \item \textbf{Attractive} for system \eqref{eq:mas_0}-\eqref{eq:mas_N} if $e(0)\notin \varepsilon(P)$ then $e(t)\to \varepsilon(P)$ as $t\to \infty$.
    \end{itemize}   
\end{Definition}
 In other words, for an invariant ellipsoid $\varepsilon(P)$, system state originates from within it will remain inside for the entire duration. Another property  is that the invariant ellipsoid is attracting if the initial position  $e(0)$ is outside the invariant set $\varepsilon(P)$ and $e(t)\to \varepsilon(P)$, as $t\to \infty$. 
 
 Since the matrix $P$ is referred to as the matrix used to characterize the ellipsoid $\varepsilon(P)$, it can be regarded as the characteristic of the impact of external disturbances on the trajectory of system. Therefore, minimize the impact of disturbance implies minimizing the invariant ellipsoid $\varepsilon(P)$. In this paper, we use the trace criterion 
 \begin{equation}
     \mathrm{tr}(P)
 \end{equation}
to characterize the size of ellipsoid, which represents the sum of the quadrates of the semiaxes of the ellipsoid invariant.

\section{Problem Statement}
 \label{sec:Pro_statement}



In this paper, we  deal with  MAS  where the communication topology is modelled by a fixed graph $\mathcal{G} = \{\mathcal{V},\mathcal{E},\mathcal{C}\}$. In $\mathcal{G}$, there is a root node which conveys its local information to the nearest nodes, and the rest nodes interact under an undirected graph. It means that $\mathcal{G}$ is composed of a root node, a set of nodes which are undirectedly connected and the directed connection from the root node to its nearest ones. In the sense of MAS, agent $0$ corresponds to the root node, which is designated as the leader, while the others, labeled from $1$ to $N$, perform as followers.  The dynamic of the leader is given in the following form 
\begin{equation}\label{eq:mas_0}
    \dot{\sigma}_0 (t) = A\sigma_0+Bu_0,\quad \sigma_0 +Bu_0 \in \R^n
\end{equation}
where the leader input $u_0$ is assumed to be known, and the followers are modelled as 
\begin{equation}\label{eq:mas_N}
    \dot{\sigma}_i (t)= A \sigma_i(t) + Bu_i(t) ,\quad i= \overline{1,N}
\end{equation}
where $\sigma_i(t)\in \R^n$ is the $i_{th}$ agent state, $u_i(t)\in\R^{m}$ is the control input of agent $i$, and the system matrix $A \in \R^{n\times n }$, $B\in \R^{ n\times m}$. In addition, the control protocol applied in this research is presented in the following form
\begin{equation} \label{eq:optimal_cont}
\resizebox{.44\textwidth}{!} 
{$
    u_i = K \sum_{j=0}^N c_{ij}(\sigma_j-\sigma_i) +u_0, \; K\in\R^{m\times n},\; i = \overline{1,N}
$}
\end{equation} 
Therefore, based on \eqref{eq:mas_0}, \eqref{eq:mas_N} and \eqref{eq:optimal_cont}, we have the error system between the state of followers and the leader as follows: 
\begin{equation}\label{eq:error_within_input_pre}
    \dot{e} = (I_N\otimes A - \tilde{\mathcal{L}}\otimes BK) e,
\end{equation}
where $e = [e_1,e_2,...,e_N]^\top$, $e_i = \sigma_i-\sigma_0$,
$$\tilde{\mathcal{L}}=\left[\begin{matrix}
l_{11} &  \ldots & l_{1N}\\
\vdots  &  \ddots  & \vdots\\
l_{N1} &  \ldots & l_{NN}\\
\end{matrix}\right]$$ 
corresponds to the undirected graph, and we let $\lambda_i\in\R$, $i=1,2,\dots,N$ be its eigenvalues. According to Lemma \ref{lem:tildeL}, $\tilde{\mathcal{L}}$ is a special case of $\hat{\mathcal{L}}$ by letting $l_{01}=\ldots=l_{0N}=0$. Thus we have $\lambda_i\geq0$. 

\begin{Definition}[inspired by \cite{olfati2004consensus}]\label{def:linearconsensus} \it
    The MAS admits the linear consensus protocol \eqref{eq:optimal_cont} if there exists a feedback gain $K$ such that the error equation \eqref{eq:error_within_input_pre} is globally asymptotically stable. 
\end{Definition} 

In this paper we study  the MAS operating in the presence of bounded external disturbances.
For example, dynamics of MAS of flying robots may be perturbed by a wind. The impact of disturbance on each agent can be considered uniformly bounded. In this paper we assume that the leader of MAS is virtual, so its dynamics is not perturbed.  
The model of such a MAS can be presented as follows: 
\begin{align}
    \dot{\sigma}_0 (t) &= A\sigma_0+Bu_0,\quad \sigma_0 +Bu_0 \in \R^n \label{eq:mas_0_dis}\\
    \dot{\sigma}_i (t) &= A \sigma_i(t) + Bu_i(t) + E\omega(t),\quad i= {1,2,\dots,N}  \label{eq:mas_N_dis}
\end{align}
where $E\in \R^{n\times p}$ is constant matrix and  $\omega(t)\in \R^p$ is the bounded external  disturbance such that 
   \begin{equation}\label{eq:bounded_dis}
       \omega^\top Q \omega \leq 1, \quad \omega \in L^{\infty}(\R_+,\R^p),\;0\prec Q=Q^{\top}\in \R^{p\times p}
   \end{equation}
Similar to the MAS without disturbances, the error system can be written as follows:
\begin{equation}\label{eq:error_with_input}
    \dot{e} =  (I_N\otimes A - \tilde{\mathcal{L}}\otimes BK) e +1_N\otimes E\omega,
\end{equation}
The main objective is to design a linear consensus protocol of MAS \eqref{eq:mas_0}, \eqref{eq:mas_N} without disturbance such that in the case of disturbance this protocol is able to  minimize the invariant ellipsoid of error system \eqref{eq:error_with_input} under constraint \eqref{eq:bounded_dis}.





\section{Main Results}
\label{Main_results}


%

\subsection{Linear Consensus Protocol in the Absence of Disturbance}




A result on linear leader-following consensus is introduced in this subsection, which is refined from \cite{ni2010leader}.


 \begin{Lemma}\it 
The following three claims are equivalent:
\begin{itemize}
    \item [a)] MAS \eqref{eq:mas_0}, \eqref{eq:mas_N} admits linear leader-following consensus,
    \item [b)] $\mathcal{G}$ has a spanning tree, and the pair $(A,B)$ is stabilizable,
    \item [c)]  Matrix Inequality 
\begin{equation}\label{LMI:linearconsensus}
\resizebox{.42\textwidth}{!} 
{$P(I_N\!\otimes\! A)\!+\!(I_N\!\otimes\! A^\top)P\!-\!P(\tilde{\mathcal{L}}\!\otimes\! BK)\!-\!(\tilde{\mathcal{L}}\!\otimes\! BK)^\top \!P\!\prec\!0$}
 \end{equation} 
 with $0\prec P\in\R^{nN\times nN}$ is feasible.
\end{itemize}
 
 
  \end{Lemma}
  
\begin{IEEEproof}
$\bullet$ a) $\Leftrightarrow$ c):
The closed-loop error equation \eqref{eq:error_within_input_pre} is globally asymptotically stable if and only if the matrix inequality \eqref{LMI:linearconsensus} is feasible. Combine Definition \ref{def:linearconsensus} about linear leader-following consensus, we can obtain that a) and c) are equivalent. 

$\bullet$ a) $\Leftrightarrow$ b): 
According to Definition \ref{def:linearconsensus}, linear leader-following consensus is considered  achieved when \eqref{eq:error_within_input_pre} is globally asymptotically stable. 
According to \cite{Fax2004TAC}, \cite{lewis2013cooperative}, the stability of 
\eqref{eq:error_within_input_pre} is equivalent to the stability of N systems $A-\lambda_i BK$, $i=1,2,\dots,N$,  where $\lambda_i$ is the eigenvalue of $\tilde{\mathcal{L}}$, $\lambda_i\geq0$. 
Inspired by Finsler's lemma \cite{boyd1994linear} and \cite{Li2012TAC}, \cite{rejeb2016synchronization}, 
if $(A, B)$ is stabilizable, then there exist solutions $(X,\gamma)\in(\R^{n\times n},\R)$ such that the following inequality holds, 
\begin{equation}\label{LMI:ctrlAB}
AX+XA^\top-\gamma BB^\top \prec0.
\end{equation}
Take $K=\frac{\gamma}{2\lambda_1}B^\top X^{-1}$, $\gamma>0$, we obtain the following inequality:
$$(A-\lambda_i BK)X+X(A-\lambda_i BK)^\top\!\!\!\prec0,\;\lambda_i\!\neq\!0,\;\forall i=1,2,\dots,N$$
Therefore, we can deduce that $\lambda_i\!>\!0$, $\forall i\!=\!1,2,\dots,N$, which is met if and only if the graph $\mathcal{G}$ has a spanning tree.
\end{IEEEproof}

\subsection{Invariant Ellipsoid of MAS with Disturbance} \label{Liear_invariant}

Consider the  leader-following system \eqref{eq:mas_0_dis}, \eqref{eq:mas_N_dis} with control protocol \eqref{eq:optimal_cont} under constrain \eqref{eq:bounded_dis},
the following theorem inspired by  \cite{nazin2007rejection} is proposed to provide a necessary and sufficient condition for being an invariant ellipsoid of system \eqref{eq:error_with_input}.

\begin{Theorem} \it 
    \label{the:lin_no_input}
    Let $P^\top\!\!\! =\!\!P\!\in\! \R^{nN\times nN}$ be positive definite.
    The ellipsoid $\varepsilon(P)$ of system \eqref{eq:error_with_input} is an invariant set under bounded external disturbance $\omega$  if and only if  there exists a real number  $\beta> 0$ such that the following matrix inequalities     holds,
\begin{equation} \label{eq:lin_LMI_op_main}
\begin{aligned}
&\resizebox{.49\textwidth}{!} 
{$ 
    \begin{bmatrix}
     P(I_N\otimes A)\!+\!(I_N\!\otimes\! A)^\top\! P  \!-\! P(\tilde{\mathcal{L}} \otimes\! BK) \!-\! \! (\tilde{\mathcal{L}}\!\otimes\! BK)^\top P\!+\! \beta P & P(1_N\!\otimes\! E) \\
     (1_N\!\otimes \!E )^\top\! P  & \!-\! \beta   Q 
\end{bmatrix}\!\!\!\preceq\! 0.
$} 
\end{aligned}
\end{equation}
\end{Theorem}
\begin{IEEEproof}
We first prove that \eqref{eq:lin_LMI_op_main} implies 
\begin{equation} \label{eq:tem_no_input_Lya}
  \begin{rcases}
  e^\top Pe \geq 1\\
  \omega^\top Q \omega \leq 1
  \end{rcases} \; \rightarrow \;
\frac{d}{dt}e^\top P e \leq 0.
\end{equation}
From \eqref{eq:lin_LMI_op_main}, one derives for any vector $[e^\top,\omega^\top]^\top$, the following inequality holds 
\begin{equation*}
\begin{bmatrix}
    e\\  \omega
\end{bmatrix}^\top \!\!\!
M
\begin{bmatrix}
    e\\  \omega
\end{bmatrix}\leq 0,
\end{equation*}
with \begin{equation*}
\resizebox{.5\textwidth}{!} 
{$ 
    M \!\!=\!\!\begin{bmatrix}
     P(I_N\otimes A)\!+\!(I_N\!\otimes\! A)^\top\! P  \!-\! P(\tilde{\mathcal{L}} \otimes\! BK) \!-\! \! (\tilde{\mathcal{L}}\!\otimes\! BK)^\top P\!+\! \beta P & P(1_N\!\otimes\! E) \\
     (1_N\!\otimes \!E )^\top\! P  & \!-\! \beta   Q 
\end{bmatrix},
$}
\end{equation*}
which implies 
\begin{align*}
  &
  \resizebox{.48\textwidth}{!} 
{$ 
e^\top\!\! (P(I_N\!\otimes\! A)\! +\!(I_N\!\otimes\! A)^\top\! \!P\!-\!P(\tilde{\mathcal{L}}\!\otimes\! BK) \!-\! (\tilde{\mathcal{L}}\!\otimes\! BK)^\top\! \!P)e
$}\\
  &\!+\!e^\top\! \!P(1_N\!\otimes\! E) \omega \! +\!\omega^\top\!\! (1_N\!\otimes\! E)^\top\!\! P e\!\leq\!  \beta(\omega^\top\!\!  Q\omega \!-\! e^\top\!\! P e).
\end{align*}
Thus for $ e^\top Pe \geq 1$,   $\omega^\top Q \omega \leq 1$ one derives  
\begin{align*}
&  \resizebox{.48\textwidth}{!} 
{$ 
     \frac{d}{dt}e^\top\!\! P e \! =\!e^\top\!\! (P(I_N\!\otimes\! A)\! +\!(I_N\!\otimes\! A)^\top\!\! P\!-\!P(\tilde{\mathcal{L}}\!\otimes \!BK) \!-\! (\tilde{\mathcal{L}}\!\otimes\! BK)^\top\!\! P)e   $}\\
     & 
\resizebox{.48\textwidth}{!} 
{$ 
+e^\top\!\! P(1_N\!\otimes\! E ) \omega \! +\!\omega^\top\!\!(1_N\!\otimes\! E)^\top\!\! P e\!\leq\! \beta(\omega^\top\! \! Q\omega \!-\! e^\top\!\! Pe) \!\leq\! 0.
$}
\end{align*}
Therefore the first claim is proved.

Next suppose that \eqref{eq:tem_no_input_Lya} holds implies the ellipsoid $\varepsilon(P)$ is not invariant, this means that there exists a trajectory starting from the initial state $e(0)$ with $g(0) = e^\top(0) P e(0) =1$ and end at an $  e(T)\in \varepsilon(P)$ such that $g(T) =e^\top(T) P e(T) >1$. Thus there exists an instant $\Bar{t}\in (0,T)$ such that  $g(\Bar{t})\geq 1$ and $\frac{d}{dt}g(\Bar{t})>0$. This contradicts the condition $\eqref{eq:tem_no_input_Lya}$. So we can conclude that $\varepsilon(P)$ is an invariant set.

Necessity:
We first prove that if $\varepsilon(P)$ is an invariant ellipsoid of system \eqref{eq:lin_LMI_op_main}, then we have 
\begin{align*}
  &
  \resizebox{.47\textwidth}{!} 
{$ 
e^\top\!\! (P(I_N\!\otimes\! A)\! +\!(I_N\!\otimes\! A)^\top\! \!P\!-\!P(\tilde{\mathcal{L}}\!\otimes\! BK) \!-\! (\tilde{\mathcal{L}}\!\otimes\! BK)^\top\! \!P)e
$}\\
  &\!+\!e^\top\! \!P(1_N\!\otimes\! E) \omega \! +\!\omega^\top\!\! (1_N\!\otimes\! E)^\top\!\! P e\! \leq\! 0,
\end{align*}
whenever $\omega^\top  Q\omega \leq e^\top P e$.
To do this, suppose this condition does not hold which means there exist $\hat{e}$ and $\hat{\omega}$ such that $\hat{\omega}^\top Q \hat{\omega} \leq \hat{e}^\top P \hat{e}$ but 
\begin{align*}
&  \resizebox{.47\textwidth}{!} 
{$ 
\hat{e}^\top\!\! (P(I_N\!\otimes\! A)\! +\!(I_N\!\otimes\! A)^\top\! \!P\!-\!P(\tilde{\mathcal{L}}\!\otimes\! BK) \!-\! (\tilde{\mathcal{L}}\!\otimes\! BK)^\top\! \!P)\hat{e}
$}\\
  &\!+\!\hat{e}^\top\! \!P(1_N\!\otimes\! E) \hat{\omega} \! +\!\hat{\omega}^\top\!\! (1_N\!\otimes\! E)^\top\!\! P \hat{e}\!>\!0.
\end{align*}
This implies $P \hat{e} \neq \textbf{0}$, thus we have $\hat{e}^\top P \hat{e} \neq 0$. Then denote $\bar{e} := \frac{\hat{e}}{\sqrt{\hat{e}^\top P \hat{e}}}$ and $\bar{\omega} := \frac{\hat{\omega}}{\sqrt{\hat{e}^\top P \hat{e}}} $.
Obviously, we have 
\begin{align*}
&  \resizebox{.47\textwidth}{!} 
{$ 
\bar{e}^\top\!\! (P(I_N\!\otimes\! A)\! +\!(I_N\!\otimes\! A)^\top\! \!P\!-\!P(\tilde{\mathcal{L}}\!\otimes\! BK) \!-\! (\tilde{\mathcal{L}}\!\otimes\! BK)^\top\! \!P)\bar{e}
$}\\
  &\!+\!\bar{e}^\top\! \!P(1_N\!\otimes\! E) \bar{\omega} \! +\!\bar{\omega}^\top\!\! (1_N\!\otimes\! E)^\top\!\! P \bar{e}\!>\!0,
\end{align*}
and $\bar{e}$ is on the boundary of $\varepsilon(P)$. For the system \eqref{eq:error_with_input} with initial state $e(0) = \bar{e}$, we derive that $g(0) :=e(0)^\top Pe(0) =1$ and $\frac{d}{dt}g(0)>0$. Therefore there must exists a small $\xi>0$, such that $g(\xi) >1$, this implies that $\varepsilon(P)$ is not an invariant set which contradicts the assumption. The first claim is proved. Thus we have 
\begin{equation*}
\resizebox{.49\textwidth}{!} 
{$ 
s^\top \!\!
\begin{bmatrix}
     P(I_N\!\otimes\! A) \!+\! (I_N\!\otimes\! A)^\top\! \!P \!-\! P(\tilde{\mathcal{L}}\!\otimes\! BK) \!-\! (\tilde{\mathcal{L}}\!\otimes\! BK)^\top\! \!P& P (1_N\!\otimes\! E) \\
     (1_N\!\otimes\! E) ^\top\!\! P  & \textbf{0}
\end{bmatrix}\!\!
s\!\leq\!  0,
$}
\end{equation*}
for $\forall s \!:\! s^\top\!\!  \begin{psmallmatrix} 0&0\\0 &Q  \end{psmallmatrix}s\!\leq\! 1$ and $s^\top \!\!\begin{psmallmatrix}
    -P &0\\0&0
\end{psmallmatrix} s \!\leq \!-1$ with $s \!=\![e^\top,\omega^\top]^\top$.
By S-procedure \cite{polyak1998convexity}, one derives 
 \eqref{eq:lin_LMI_op_main} holds for $\beta\geq 0$. For the case $\beta =0$, the following inequality must hold, 
\begin{equation} 
\begin{aligned}
&\resizebox{.49\textwidth}{!} 
{$ 
    \begin{bmatrix}
     P(I_N\!\otimes\! A)\!+\!(I_N\!\otimes\! A)^\top\!\! P  \!-\! P(\tilde{\mathcal{L}} \!\otimes\! BK) \!-\!  (\tilde{\mathcal{L}}\!\otimes\! BK)^\top\!\!P  & P(1_N\!\otimes \!E) \\
     (1_N\!\otimes\! E )^\top\!\!P   &\textbf{0} 
\end{bmatrix}\!\!\preceq\! 0.
$} 
\end{aligned}
\end{equation}
The latter inequality  holds only if $E=\textbf{\zero}$, which contradicts the assumption of this theorem. Therefore  we obtain that $\beta $ should be positive.
\end{IEEEproof}

\begin{Corollary}\label{cor:lin_attractive} \it
Let the gain $K$, as defined in \eqref{eq:optimal_cont}, be selected such that the MAS \eqref{eq:mas_0}, \eqref{eq:mas_N} achieves leader-following consensus without  disturbances. Then any invariant ellipsoid of perturbed  leader-following system \eqref{eq:error_with_input}    is attractive.
\end{Corollary}
\begin{IEEEproof}
Since the leader-following consensus of system    is achieved, then one derives for the system \eqref{eq:error_within_input_pre} without disturbance
\begin{align*}
  \resizebox{.49\textwidth}{!} 
{$ 
  \frac{d}{dt}e^\top\!\! P e \! =\! 
e^\top\!\! (P(I_N\!\otimes\! A) \!+\!(I_N\!\otimes\! A)^\top\!\! P \!-\! P(\tilde{\mathcal{L}}\!\otimes \!BK) \!-\! (\tilde{\mathcal{L}}\!\otimes\! BK)^\top\! \!P)e \!\leq\! 0.
$}
\end{align*}
From Theorem \ref{the:lin_no_input}, any invariant ellipsoid of perturbed leader-following system \eqref{eq:error_with_input} implies that there exists a  $\beta >0$ such that 
\begin{align*}
  &
  \resizebox{.48\textwidth}{!} 
{$ 
e^\top\!\! (P(I_N\!\otimes\! A)\! +\!(I_N\!\otimes\! A)^\top\! \!P\!-\!P(\tilde{\mathcal{L}}\!\otimes\! BK) \!-\! (\tilde{\mathcal{L}}\!\otimes\! BK)^\top\! \!P)e
$}\\
  &\!+\!e^\top\! \!P(1_N\!\otimes\! E) \omega \! +\!\omega^\top\!\! (1_N\!\otimes\! E)^\top\!\! P e\!\leq\!  \beta(\omega^\top\!\!  Q\omega \!-\! e^\top\!\! P e).
\end{align*}
Therefore, we obtain that for the system \eqref{eq:error_with_input} with disturbance
\begin{align*}
&  \resizebox{.48\textwidth}{!} 
{$ 
     \frac{d}{dt}e^\top\!\! P e \! =\!e^\top\!\! (P(I_N\!\otimes\! A)\! +\!(I_N\!\otimes\! A)^\top\!\! P\!-\!P(\tilde{\mathcal{L}}\!\otimes \!BK) \!-\! (\tilde{\mathcal{L}}\!\otimes\! BK)^\top\!\! P)e   $}\\
     & 
\resizebox{.48\textwidth}{!} 
{$ 
+e^\top\!\! P(1_N\!\otimes\! E ) \omega \! +\!\omega^\top\!\!(1_N\!\otimes\! E)^\top\!\! P e\!\leq\! \beta(\omega^\top\! \! Q\omega \!-\! e^\top\!\! Pe) \!<\! 0
$}
\end{align*}
holds whenever $e^\top P e> 1$ and $\omega^\top  Q\omega\leq 1$, which is 
sufficient condition of being an attractive ellipsoid of system \eqref{eq:error_with_input}. 
\end{IEEEproof}

From the practical application point of view, it is natural to introduce the constraint of input within the framework of the proposed approach. 
\begin{Corollary} \it 
    The linear control protocol satisfies the following constraint 
    \begin{equation}
        (u(e)-u_0)^\top (u(e)-u_0) \leq \eta^2, \quad \forall e \in \R^{nN}  :\quad e^\top P e\leq 1,
    \end{equation}
    for some $\eta >0$ if and only if the following LMI  holds,
    \begin{equation}\label{eq:Lin_input_constraint_LMI}
        \begin{bmatrix}
             P & \tilde{\mathcal{L}}^\top \!\!\otimes\! K^\top\!\\
            \tilde{\mathcal{L}}\otimes K  & \eta^2 I_{nN}
        \end{bmatrix} \succeq 0.
    \end{equation}
\end{Corollary}
\begin{IEEEproof}
From \eqref{eq:optimal_cont}, one derives  that $u(e) = (\tilde{\mathcal{L}}\otimes K)e$ \cite{ni2010leader}. To satisfy the control constraints, it requires
\begin{equation}
    e^\top (\tilde{\mathcal{L}}\otimes K)^\top (\tilde{\mathcal{L}}\otimes K) e \leq  \eta^2,\quad \forall e :\quad e^\top P e\leq 1,
\end{equation}
which is a classical S-procedure for two quadratic forms. It is equivalent to there exists $\tau \geq 0$ such that 
\begin{equation}
   (\tilde{\mathcal{L}}\otimes K)^\top (\tilde{\mathcal{L}}\otimes K) \preceq \tau P,\quad \tau \leq \eta^2.
\end{equation} 
Since the minimal invariant ellipsoid is interesting to us, we take 
\begin{equation}
    \tau =\tau_{max} = \eta^2,
\end{equation}
then one derives
\begin{align}
& (\tilde{\mathcal{L}}\otimes K)^\top (\tilde{\mathcal{L}}\otimes K) \preceq \eta^2 P.
\end{align}
Use Schur complement, we obtain \eqref{eq:Lin_input_constraint_LMI}.  \end{IEEEproof}

\begin{Corollary} \label{lem:worst_case} \it
The worst case of external disturbance $\omega^\ast $ for system \eqref{eq:error_with_input} can be given in the following form 
\begin{equation} \label{eq:worst_dis}
\omega^\ast = \frac{Q^{-   1}(1_N\otimes E)^\top P e}{\|Q^{-   \frac{1}{2}}(1_N\otimes E)^\top P e\|}.
\end{equation}
\end{Corollary}
\begin{IEEEproof}
    The worst disturbance will push the system trajectories to the boundary of the invariant ellipsoid, which requires 
    \begin{equation}
      \frac{d}{dt} e^\top P e \quad \to \quad \mathrm{Max}.
    \end{equation}
Since 
\begin{align*}
&\resizebox{.48\textwidth}{!} 
{$
\frac{d}{dt} e^\top\!\!P e \!=\!    
e^\top \!\!(P(I_N\!\otimes\! A) \!+\!(I_N\!\otimes\! A)^\top\! \!P \!-\! P(\tilde{\mathcal{L}}\!\otimes\! BK) \!-\! (\tilde{\mathcal{L}}\!\otimes\! BK)^\top\! \!P)e 
$}\\
&\!+\!e^\top\!\! P(1_N\!\otimes \!E )\omega \! +\!\omega^\top\! \!(1_N\!\otimes\! E)^\top\!\! P e.    
\end{align*}
Thus it becomes the problem 
\begin{equation}\label{eq:pro_worste}
  \max_{\omega^\top Q \omega =1} \langle\omega,(1_N\otimes E)^\top P e\rangle.
\end{equation}
Denote $\Bar{\omega} = Q^{\frac{1}{2}}\omega$, since $0\prec Q=Q^\top$, \eqref{eq:pro_worste} can be rewritten as 
\begin{equation}
    \max_{\Bar{\omega}^\top \Bar{\omega}=1}\langle\Bar{\omega},Q^{-   \frac{1}{2}}(1_N\otimes E)^\top P e\rangle.
\end{equation}
The obvious solution of $\Bar{\omega}$ is 
\begin{equation}
    \Bar{\omega}^\ast = \frac{Q^{-   \frac{1}{2}}(1_N\otimes E)^\top P e}{\|Q^{-   \frac{1}{2}}(1_N\otimes E)^\top P e\|}.
\end{equation}
Finally one derives the solution $\omega^\ast$ as follows
\begin{equation}
    \omega^\ast  = \frac{Q^{-   1}(1_N\otimes E)^\top P e}{\|Q^{-   \frac{1}{2}}(1_N\otimes E)^\top P e\|}.
\end{equation}
\end{IEEEproof}

\subsection{Minimization of invariant ellipsoid}

Since \eqref{eq:lin_LMI_op_main}  guarantees $\varepsilon(P)$ to be attractive/invariant ellipsoid, and the trace of matrix $P$ characterizes the size of invariant ellipsoid (represent the impact of external disturbance). Therefore, the external disturbance of system \eqref{eq:error_with_input} can be optimally rejected by solving the following optimal problem: 
 \begin{equation}\label{eq:min}
 \begin{aligned}
  &\mathrm{tr}(X) \quad \to \quad \mathrm{Min} \quad \text{where } X=P^{-1} \\
  &\text{under constraints  \eqref{eq:lin_LMI_op_main} and \eqref{eq:Lin_input_constraint_LMI}},
  \end{aligned}
 \end{equation}
 where the constraint  \eqref{eq:Lin_input_constraint_LMI} is added to avoid the infinite input. 
The following corollary allows one to limit the search for the minimum invariant ellipsoid to the one-parameter family \eqref{eq:lin_LMI_op_main}, simplifying the problem to a one-dimensional  minimization over a finite interval.
\begin{Corollary} \it 
For any given matrix gain $K$  that makes $ (I_N\otimes A-\Tilde{\mathcal{L}}\otimes BK)$ Hurwitz, the minimal (in the sense of \eqref{eq:min}) invariant ellipsoid of system \eqref{eq:lin_LMI_op_main} belongs to the one parameter family of ellipsoid generated by matrix $P$ and $K$ that satisfying  
\begin{equation}\label{eq:convex_eq}
\begin{aligned}
&\resizebox{.49\textwidth}{!} 
{$ 
    \begin{bmatrix}
     P(I_N\!\otimes\! A)\!+\!(I_N\!\otimes\! A)^\top\!\! P \!-\!P(\tilde{\mathcal{L}}\! \otimes\! BK)\! -\!  (\tilde{\mathcal{L}}\!\otimes\! BK)^\top\!\! P\!+\! \beta P& P(1_N\!\otimes\! E) \\
     (1_N\!\otimes \!E )^\top\!\! P  & - \beta Q 
\end{bmatrix}\!\!=\!0,
$} 
\end{aligned}
\end{equation}
over the interval $0<\beta <-2\mathrm{Max}(\mathrm{Re}\lambda_i(I_N\otimes A-\Tilde{\mathcal{L}}\otimes BK))$. 
\end{Corollary}
\begin{IEEEproof}
  Follow the Lemma A.16 of \cite{polyak2002robastnaya}: Let $A$  be Hurwitz matrix and $(A, E)$ be controllable, then  the solution of Lyapunov equation 
\begin{equation}
    AX+XA^\top +EE^\top =0
\end{equation}
is the solution of optimal problem
\begin{equation}
    tr(X)\quad\to\quad\mathrm{Min},
\end{equation}
under constraint 
\begin{equation}
    AX+XA^\top +EE^\top \preceq 0.
\end{equation}

Then the equation 
\eqref{eq:convex_eq} can be rewritten as 
\begin{equation}\label{eq:convex_eq_P}
\begin{aligned}
&\resizebox{.48\textwidth}{!} 
{$
P(I_N\!\otimes\! A\!-\!\Tilde{\mathcal{L}}\!\otimes\! BK\!+\! \frac{\beta }{2}I_{nN} ) \! +\! (I_N\!\otimes\! A\!-\!\Tilde{\mathcal{L}}\!\otimes\! BK \!+\! \frac{\beta }{2}I_{nN} )^\top\!\! P $} \\
&  \!+\!\frac{1}{\beta}P(1_N\!\otimes\! E)Q(1_N\!\otimes\! E)^\top\! \!P\!    =\!0.
\end{aligned}
\end{equation}
Similarly, we require that 
$(I_N\otimes A-\Tilde{\mathcal{L}}\otimes BK+ \frac{\beta }{2} I_{nN})$ is Hurwitz matrix 
with
\begin{equation}
\mathrm{Re}\lambda_i(I_N\otimes A-\Tilde{\mathcal{L}}\otimes BK+ \frac{\beta }{2}I_{nN} )<0,
\end{equation}
then we obtain that \eqref{eq:convex_eq} has unique solution when $0 <\beta<-2\mathrm{Max}(\mathrm{Re}\lambda_i(I_N\otimes A-\Tilde{\mathcal{L}}\otimes BK))$.

\end{IEEEproof}

\section{Simulation Examples}

\label{sim_results}
\textcolor{blue}{
In the following simulation example, we firstly define a MAS model with bounded external disturbance  and an input constraint. Secondly, the optimal problem \eqref{eq:min} is solved to obtain the optimal solution $K$. Finally, the MAS is simulated based on the optimal linear control protocol \eqref{eq:optimal_cont}.    
}
\subsection{Numerical Results}
\begin{figure}[ht]    
\centering
\includegraphics[scale=1]{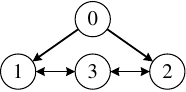}  
\caption{ Topology of MAS  }
\label{fig:MAS_example_1}
\end{figure}

\begin{figure}[ht]    
\centering
\includegraphics[scale=0.25]{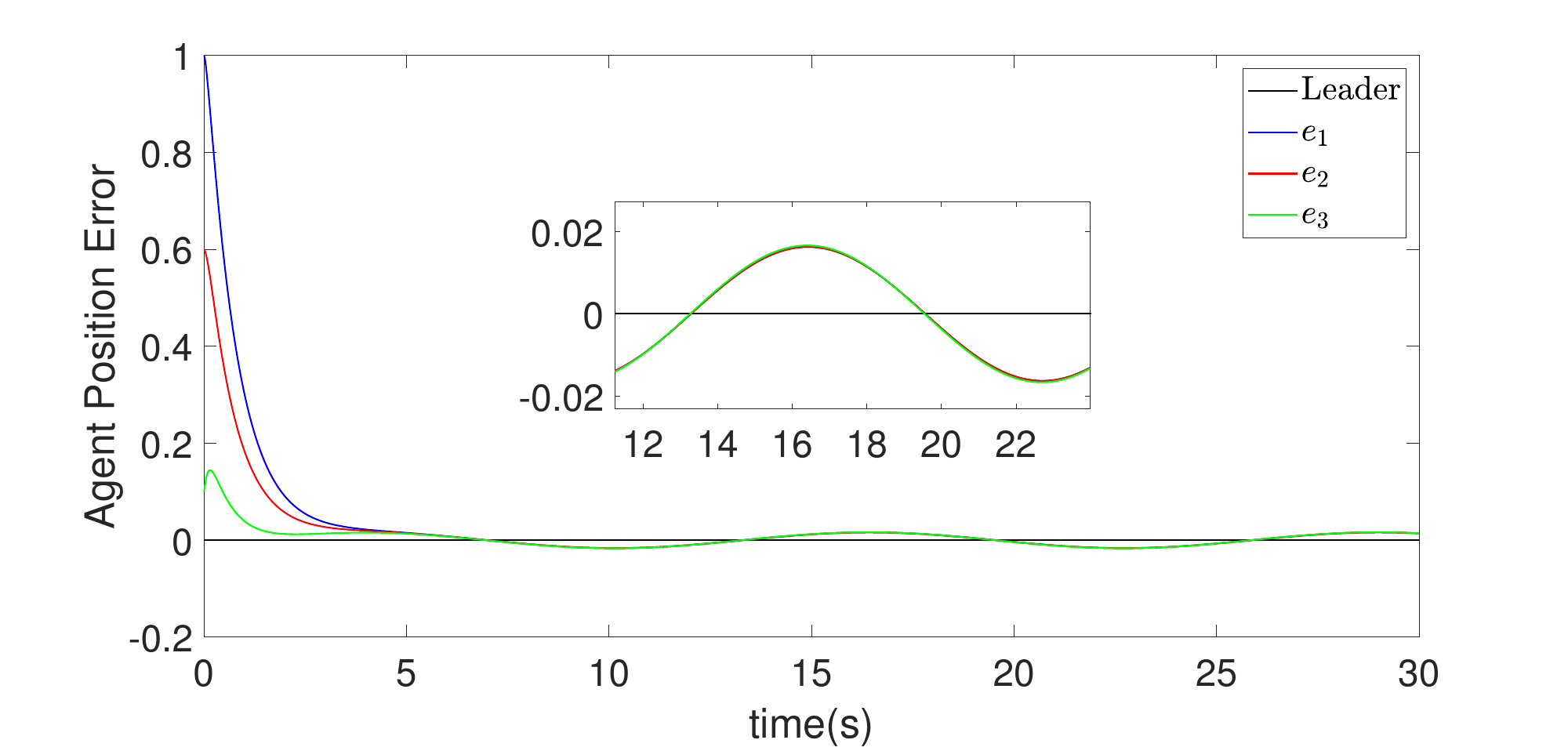}  
\caption{ Consensus trajectory of three agents under bounded disturbance  }
\label{fig:position_tra_1}
\end{figure}

Consider the MAS \eqref{eq:mas_0}, \eqref{eq:mas_N} with 4 agents regulated by the
linear consensus protocol \eqref{eq:optimal_cont}. The  fixed communication topology is depicted in Fig.\ref{fig:MAS_example_1}. 
The agents are labeled by $i = 0,1,2,3$, where the agent $0$ is the leader and
the agent $1,2,3$ are followers. In this example,  the dynamic of 
leader and follower are modeled by   \eqref{eq:mas_0} and \eqref{eq:mas_N} respectively. The dynamic system
matrix A,  B and E are randomly chosen as follows:
\begin{align} \label{eq:system_sim}
A = \begin{bmatrix}
    0 & 1\\
    0 & 0
\end{bmatrix},\quad  B = \begin{bmatrix}
    0\\ 1
\end{bmatrix},\quad E =
\begin{bmatrix}
  1 &0  \\
 0& 1  
\end{bmatrix}
\end{align}
Let the weight $\omega_{ij} = 1$ if the agent $i$ receives the information from agent $j$, otherwise $\omega_{ij} = 0$. Thus the Laplacian matrix can be written as follows
\begin{equation}
  \mathcal{L} = \begin{bmatrix}
      0 & 0 & 0 & 0\\
      -1 & 2 & 0 & -1\\
      -1 & 0 & 2 & -1 \\
      0 & -1 & -1& 2
  \end{bmatrix}. 
\end{equation}
Then one derives 
\begin{equation}
    \Tilde{\mathcal{L}} = 
    \begin{bmatrix}
       2 & 0 & -1\\
       0 & 2 & -1 \\
       -1 & -1& 2
  \end{bmatrix}. 
\end{equation}
The external disturbance is bounded with $\omega^\top Q \omega \leq 1$ where 
\begin{equation}
    Q = \begin{bmatrix}
        800    &0\\
        0&   4000
    \end{bmatrix},\quad
    \omega = \begin{bmatrix}
       \frac{1}{50}\\ \frac{1}{80} 
    \end{bmatrix}\sin{\frac{t}{2}}.
\end{equation}
The bounded input is defined by $\eta = 50000$.  
Solve the optimal problem \eqref{eq:min} with bounded input \eqref{eq:Lin_input_constraint_LMI} by the solver {\it BMIBNB } of Yalmip \cite{yalmip_BMI}, we obtain the optimal control protocol 
\begin{align}
&K = \begin{bmatrix}
       46.6001  & 25.6217
    \end{bmatrix},\\
   & P = 10^{3}\!\times\!\begin{bsmallmatrix}
    1.9963  &  0.0008 &  -1.2544 &  -0.0003 &  -0.6919 &  -0.0018 \\
    0.0008  &  0.0188 &  -0.0005 & -0.0135  &  0.0014  & -0.0037\\
   -1.2544  & -0.0005 &   1.9291 &  -0.0000 &  -0.6245 &  -0.0027\\
   -0.0003  & -0.0135 &  -0.0000 &   0.0186 &   0.0030 &  -0.0030\\
   -0.6919  &  0.0014 &  -0.6245 &   0.0030 &   1.2549 &   0.0049\\
   -0.0018  & -0.0037 &  -0.0027 &  -0.0030 &   0.0049 &   0.0080  
    \end{bsmallmatrix}.    
\end{align}
{\color{blue} In the following simulation examples, the optimal $K$ above is applied.}  
In the first example, the leader is assumed to be at origin. The initial state of MAS is given by 
$$\sigma = [\sigma_0,\sigma_1,\sigma_2,\sigma_3]^\top =\begin{bmatrix}
    0, 0,1,0,0.6,0,0.1,0.5
\end{bmatrix}^\top,$$
with $u_0 ={0}$. Fig. \ref{fig:position_tra_1} shows the consensus trajectory of three different agents under bounded disturbance.
It is obvious to see that the trajectory converges to certain invariant set centered by the leader. \textcolor{blue}{The maximum position error between leader and followers is from the agent 3, the maximum absolute error is about $0.01659$. This means the position errors of followers converges to the invariant set $[-0.01659,0.01659]$.}

\begin{figure}[ht]    
\centering
\includegraphics[scale=0.25]{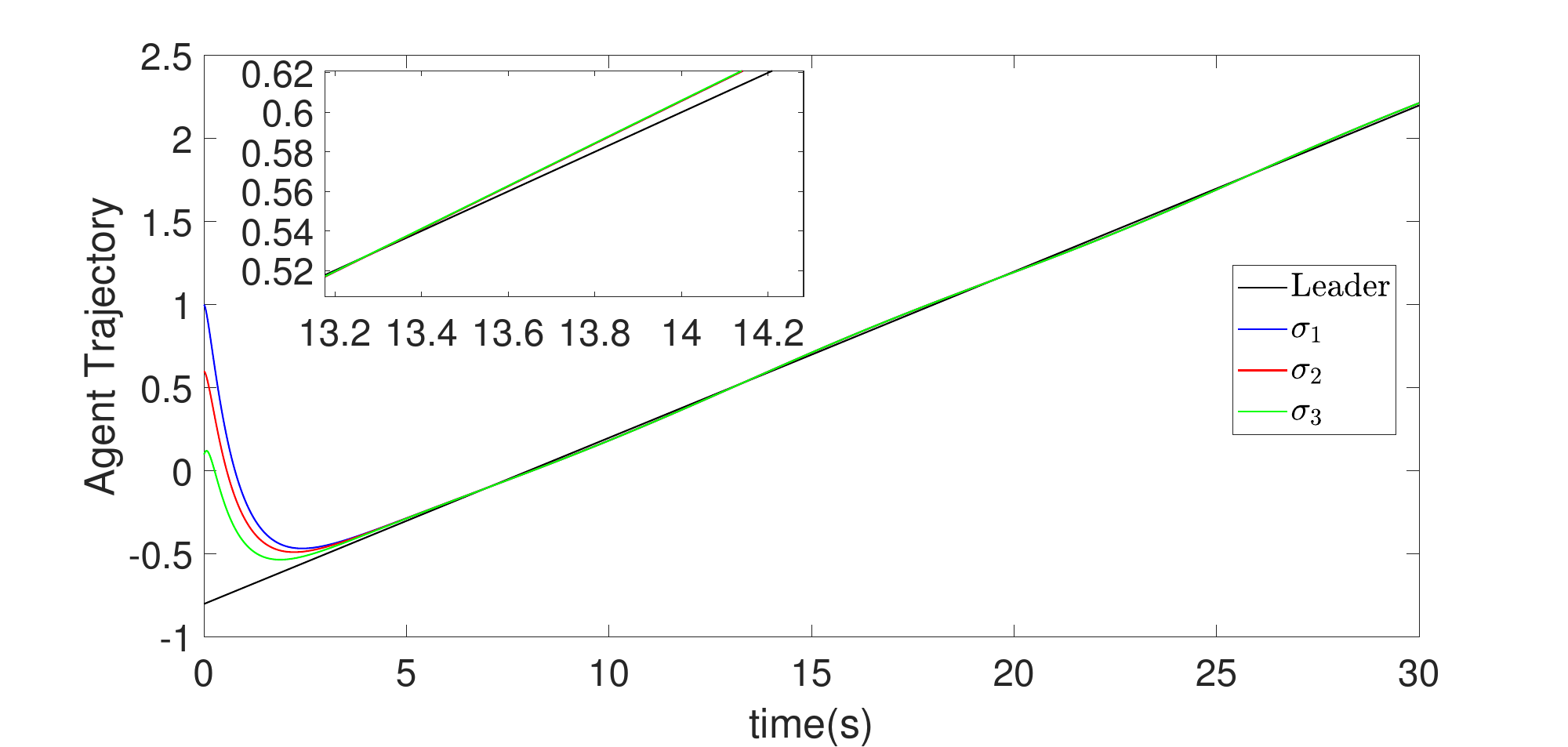}  
\caption{Trajectory of three agents tracking leader under bounded disturbance  }
\label{fig:position_tra_2}
\end{figure}

In the second example, it is assumed that the leader is moving over time. The initial state of MAS is given by 
\begin{equation*}
\resizebox{.48\textwidth}{!} 
{$\sigma = [\sigma_0,\sigma_1,\sigma_2,\sigma_3]^\top =\begin{bmatrix}
    -0.8,0.1,1,0,0.6,0,0.1,0.5
\end{bmatrix}^\top, 
$}
\end{equation*}
with $u_0 = 0.01$. Figure \ref{fig:position_tra_2} depicts the trajectories of three agents as they follow the leader. {\color{blue}
Similarly to the previous example, the maximum absolute position error comes from the agent 3, which is about $0.0243$.}

In the third example, the MAS works in the case of worst disturbance as defined  in \eqref{eq:worst_dis}. To make the trajectory  more clear, the following initial state is applied 
$$\sigma = [\sigma_0,\sigma_1,\sigma_2,\sigma_3]^\top =\begin{bmatrix}
    0,0,1,0,0.6,0,0.1,0.5
\end{bmatrix}^\top,
$$
with $u_0 =0$.
  
\begin{figure}[ht]    
\centering
\includegraphics[scale=0.25]{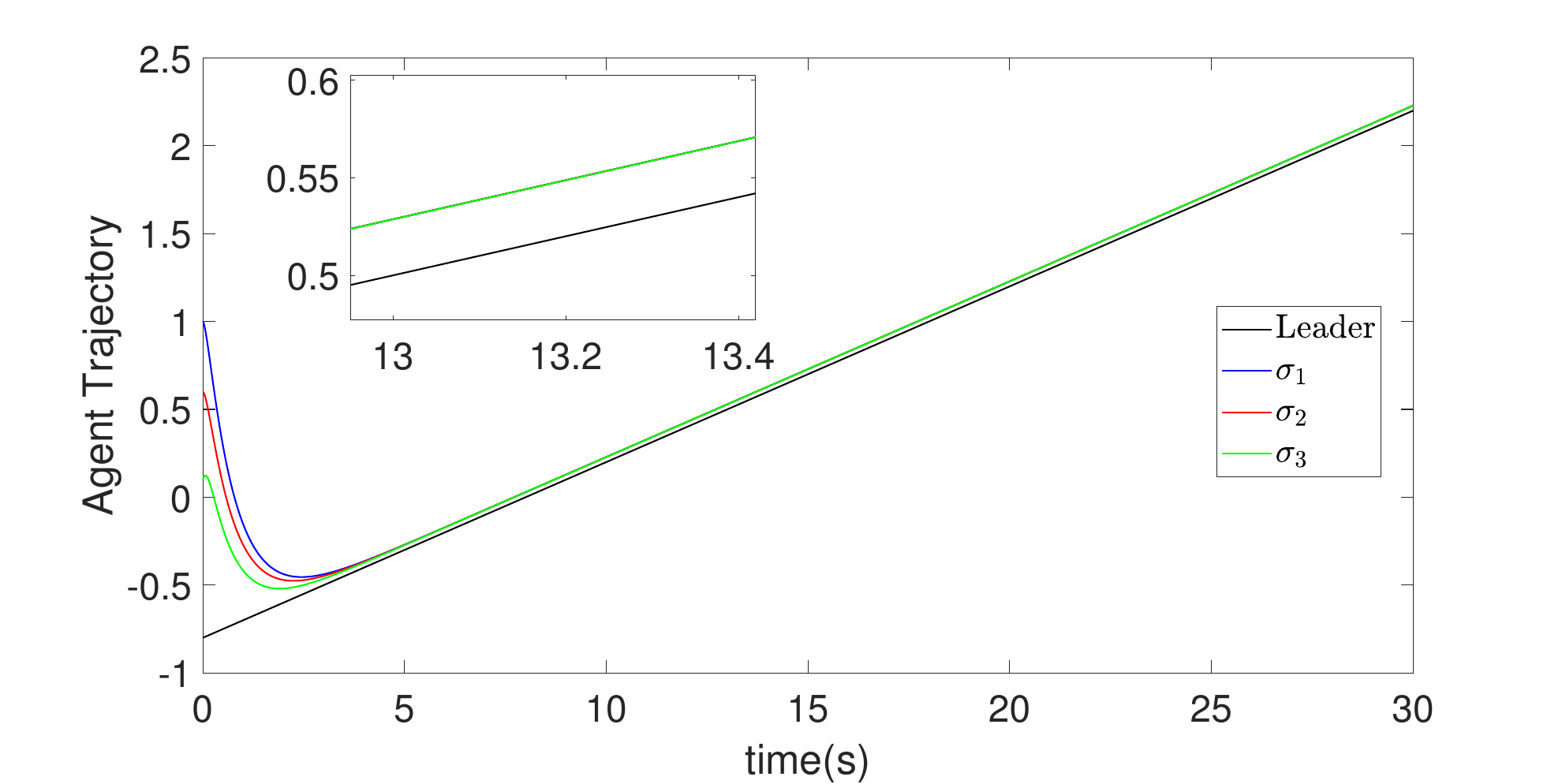}  
\caption{Trajectory of three agents tracking leader under worst disturbance  }
\label{fig:position_tra_3}
\end{figure}

It is evident from Figure \ref{fig:position_tra_3} that the worst disturbance pushes the system's trajectory farthest away from the reference position. {\color{blue} This maximum absolute position error is constant about $0.0397$.}

\section{Conclusion}\label{sec:conclusion}
This article extends the invariant/attractive ellipsoid method \cite{nazin2007rejection,khlebnikov2011optimization,Poznyak_etal2014:Book} to the MAS.  The optimal control protocol is derived through solving semidefinite matrix inequalities, which characterize the minimal invariant/attractive ellipsoid. From the practical application point of view, a bounded input is introduced. Simulation results demonstrate that the MAS state converges to a minimal invariant ellipsoid in the presence of bounded disturbances. Remarkably, in the worst case of disturbance, the system trajectories converge and then remain inside the invariant ellipsoid.  
This outcome signifies the effectiveness of the optimal control protocol in minimizing the impact of disturbances on the MAS. 

\bibliographystyle{IEEEtran}
\bibliography{Invariant_set_Ref}\ 

\end{document}